\documentclass[11pt]{amsart}
\usepackage[T1]{fontenc}
\usepackage[utf8]{inputenc}
\usepackage[english]{babel}
\usepackage{csquotes} 
\usepackage[a4paper,margin=3cm]{geometry}
\usepackage{graphicx}
\usepackage[backend=biber,style=numeric,bibencoding=utf8,language=auto,autolang=other,giveninits=true,doi=false,isbn=false,url=false,maxnames=10,sorting=nty]{biblatex}
\usepackage{caption} 
\usepackage[hyperfootnotes=false]{hyperref} 
\hypersetup{
	colorlinks = true,
	linkcolor = {blue},
	urlcolor = {red},
	citecolor = {blue}
}
\usepackage{footmisc}
\usepackage[nameinlink,capitalise,sort,noabbrev]{cleveref}
\crefname{equation}{}{} 
\crefname{enumi}{}{} 

\usepackage{amsmath}
\usepackage{amsthm}
\usepackage{amssymb}
\usepackage{esint} 
\usepackage{mathrsfs} 
\usepackage{faktor} 
\usepackage{dsfont} 
\usepackage[dvipsnames]{xcolor}
\usepackage{array}
\usepackage{hhline}
\usepackage{enumitem} 
\usepackage{comment} 
\usepackage[toc,page]{appendix} 
\usepackage{imakeidx} 
\usepackage{xparse} 
\usepackage{mathtools}
\usepackage{fancyhdr} 
\usepackage{ifthen} 
\usepackage{forloop} 
\usepackage{xstring}
\usepackage{tikz}
\usetikzlibrary{babel} 
\usetikzlibrary{cd} 
\usepackage{emptypage} 
\usepackage{listings} 
\usepackage{mathabx} 

\theoremstyle{plain}
\newtheorem{lemma}{Lemma}[section]
\newtheorem{proposition}[lemma]{Proposition}
\newtheorem{theorem}[lemma]{Theorem}
\newtheorem{corollary}[lemma]{Corollary}

\theoremstyle{definition}

\theoremstyle{remark}
\newtheorem{remark}[lemma]{Remark}

\numberwithin{equation}{section}

\renewcommand{\d}{\mathrm{d}}

\newcommand{\R}{\mathbb{R}}

\newcommand{\assign}{:=}

\newcommand{\tmop}[1]{\ensuremath{\operatorname{#1}}}

\begin{document}

\title[]{Nested discontinuous asymptotic profiles for the viscous Burgers equation with infinite mass}



\author[N.~De Nitti]{Nicola De Nitti}
\address[N.~De Nitti]{EPFL, Institut de Mathématiques, Station 8, 1015 Lausanne, Switzerland.}
\email{nicola.denitti@epfl.ch}

\author[E.~Pacherie]{Eliot Pacherie}
\address[E.~Pacherie]{CNRS and CY Cergy Paris Universit\'e , Laboratoire AGM, Avenue Adolphe Chauvin 2, 95302 Cergy-Pontoise, France.}
\email{eliot.pacherie@cyu.fr}

\begin{abstract}
We study the viscous Burgers equation with a family of initial data having infinite mass. After rescaling, the solution converges toward a bounded discontinuous profile in the long-time limit. Moreover, by changing the scale near the discontinuity point again, we find a new profile that is also discontinuous. This process can be repeated an arbitrary number of times. 
\end{abstract}

\maketitle

\section{Introduction}
\label{sec:intro}

We consider the Cauchy problem associated with the  \emph{viscous Burgers equation}, which was first introduced by H.~Bateman \cite{B} and J.~M.~Burgers \cite{MR1147, MR27195} as an approximation for equations of fluid-flow: 
\begin{align}\label{eq:burgers}
\begin{cases}
\partial_t f + f
\partial_x f - \partial_x^2 f = 0, & (x,t) \in \R \times (0,+\infty), \\  
f(x,0) = f_0(x), & x \in \R, 
\end{cases}
\end{align} 
with initial condition satisfying 
\begin{align}\label{eq:ic}
    \begin{aligned}
 &  f_0 \in
C^1_{\tmop{loc}} (\mathbb{R}; \mathbb{R}), \\
& f_0 (x) \sim \frac{\kappa}{| x |^{\alpha}}, \qquad f_0' (x) \sim \frac{- \alpha
   \kappa}{| x |^{1 + \alpha}}, \qquad \text{ as $ x \to \pm \infty$,} \\
& \text{for some $\kappa > 0$ and $\alpha \in \left] 0, 1\right[ $.}
    \end{aligned}
\end{align} 

While this function is not integrable, using the Hopf--Cole formula \cite{MR42889,MR47234} we can show that the solution of \cref{eq:burgers} with this initial datum is globally well-posed for positive times. 
 
We are interested in studying the long-time behavior of the solution $f$ of this Cauchy problem. 

The problem of the long-time asymptotics for advection-diffusion equations with has received much attention. When the initial datum is integrable and of mass $M$, the solutions to the equations under consideration often resemble the associated self-similar profile of mass $M$; however, they may also exhibit other types of asymptotic behavior (e.g., weakly non-linear, linear, or strongly non-linear) depending on the form of the convective term. We refer, e.g., to \cite{zuazua2020asymptotic,MR1926919,MR1124296,MR1266100,MR1233647,MR1443609,MR1440033,MR4725026,MR4659287,MR2346798} for some results on this topic.

When the initial data is not integrable, but merely bounded, more complicated asymptotic behavior may arise (cf., e.g., \cite{zbMATH01799645,ThierryGallay2023}). The source of difficulty and interest in our study is, in fact, the lack of integrability and the particular structure of \cref{eq:ic}.

As $t \nearrow +\infty$, the solution of \cref{eq:burgers} converges to $f_\infty \equiv 0$ (in $L^{\infty} (\mathbb{R})$, for instance). Upon rescaling, we can see the
asymptotic profile: as shown in {\cite[Theorem 1.5]{MR4687278}}, there exists $z_c \in
\mathbb{R}$ depending on $\alpha, \kappa$ such that
\[ t^{\frac{\alpha}{1 + \alpha}} f \left( t^{\frac{1}{1 + \alpha}} z, t
   \right) \longrightarrow \mathfrak{p} (z), \quad \text{for any $z \neq z_c$ as $t \nearrow + \infty$,} \]
where $\mathfrak{p}$ is a
bounded function. The convergence is uniform in $z$ as
long as we remove a neighborhood of $z_c$. The function $\mathfrak{p}$ is
smooth except at $z_c$, where it is discontinuous. The discontinuity of $\mathfrak{p}$ at $z_c$ can be seen as a form of boundary
layer, separating the area $z < z_c$ where the effect of the tail of $f_0$ at
$- \infty$ dominated the effect of the tail at $+ \infty$, and the area $z >
z_c$ where the roles are reversed. 

\begin{figure}[ht]
    \centering
    \includegraphics[width=9cm]{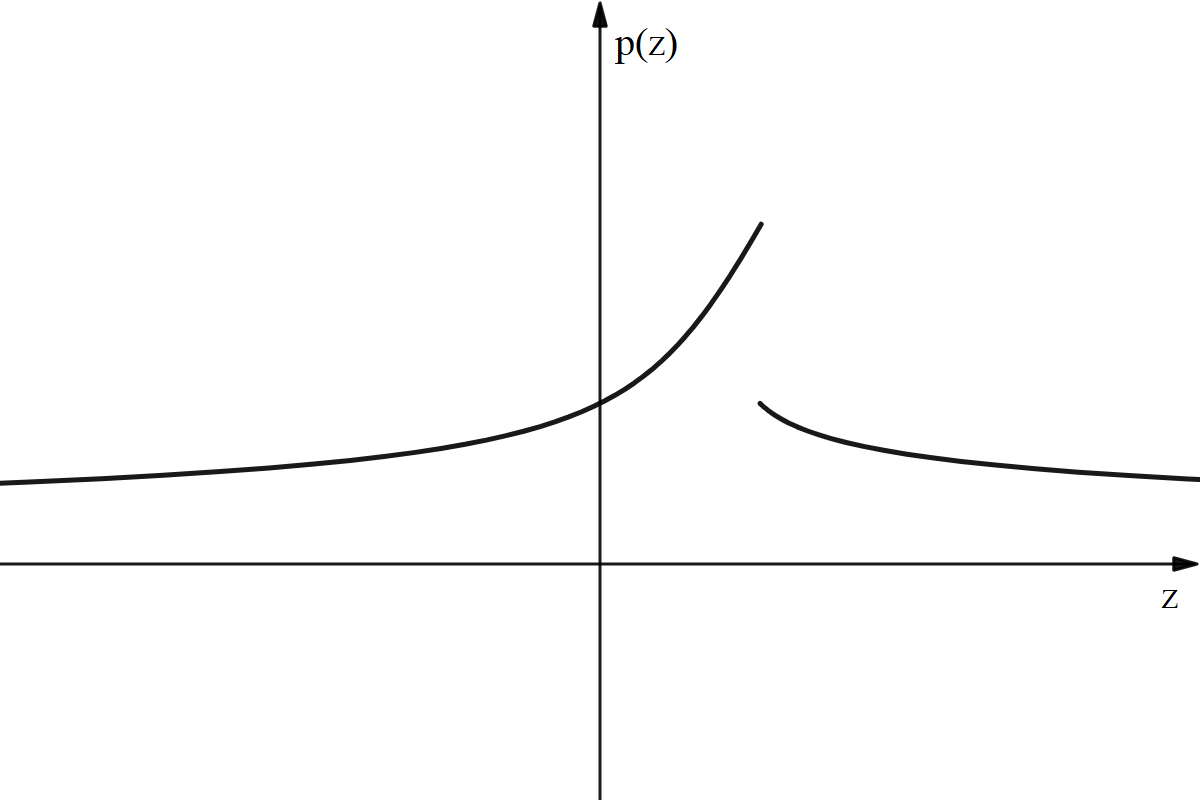}
    \caption{Graph of the function $\mathfrak{p}$. The point of discontinuity is $z_c$.}
    \label{fig:profile}
\end{figure}

It has been shown in {\cite[Theorem 1.4]{ghoul2024nonlinear}} that, for some specific initial data $f_0$, 
\[ 
\begin{aligned} &t^{\frac{\alpha}{1 + \alpha}} f \left( t^{\frac{1}{1 + \alpha}} \left( z_c
   + t^{- \frac{1 - \alpha}{1 + \alpha}} x \right), t \right) \longrightarrow
   \mathfrak{a} \tanh (2\mathfrak{a}x) + \frac{\mathfrak{p} (z_c^+)
   +\mathfrak{p} (z_c^-)}{2}, \\ &\text{for all $x \in \mathbb{R}$ when $t \rightarrow + \infty$,} 
\end{aligned} \]
   where 
\[ \mathfrak{a} \assign \frac{\mathfrak{p} (z_c^+) -\mathfrak{p} (z_c^-)}{2}
\]
As $x
\rightarrow \pm \infty$, this profile converges to $\mathfrak{p}
(z_c^{\pm})$, respectively. In other words, we can do a rescaling centered at the discontinuity of $\mathfrak{p}$ in such a way to see a new, smooth, and bounded profile. The limits of this profile at $\pm \infty$ are $\mathfrak{p}(z_c^{\pm})$, the values on both sides of the discontinuity.  We note that, in the particular case of \cite[Theorem 1.4]{ghoul2024nonlinear}, $f_0 (x) - \frac{\kappa}{| x |^{\alpha}}$ is integrable when $x \rightarrow \pm \infty$.

In this paper, our goal is to show that this is not a generic behavior. Let us start with the following example.

\begin{proposition}
  \label{prop:ab} 
  Let us consider the solution of the viscous Burgers equation
 \cref{eq:burgers} with initial data 
 \begin{align}\label{eq:ic1}
    \begin{aligned}
 &  f_0 \in
C^1_{\tmop{loc}} (\mathbb{R}; \mathbb{R}), \\
&f_0 (x) = \frac{\kappa_1}{| x |^{\alpha}} + \frac{\kappa_2}{| x
     |^{\beta}} + o \left( \frac{1}{| x |^{\beta}} \right) \quad \text{ as $x \rightarrow \pm \infty$},  \\  &f_0' (x) = -
     \frac{\alpha \kappa_1}{| x |^{1 + \alpha}} - \frac{\beta \kappa_2}{| x
     |^{1 + \beta}} + o \left( \frac{1}{| x |^{1 + \beta}} \right) \quad \text{ as $x \rightarrow \pm \infty$}, \\
    & \text{for some  $\kappa_1, \kappa_2 > 0$, $\alpha \in]
  0, 1 [$, and $\beta \in \left] \alpha, \frac{1 + \alpha}{2} \right[$}.
  \end{aligned}
     \end{align}
Then,  for any $x \neq 0$,
  \[ t^{\frac{\beta - \alpha}{1 + \alpha}} \left( t^{\frac{\alpha}{1 +
     \alpha}} f \left( t^{\frac{1}{1 + \alpha}} \left( z_c + t^{- \frac{\beta
     - \alpha}{1 + \alpha}} x \right), t \right) -\mathfrak{p} \left( z_c +
     t^{\frac{\beta - \alpha}{1 + \alpha}} x \right) \right) \longrightarrow
     \mathfrak{q} (x) \quad \text{as $t \nearrow + \infty$}, \]
where $\mathfrak{q}$ is constant on $] -  \infty, 0 [$ and on $] 0, + \infty [$, and is discontinuous at $x = 0$ for
  almost every $\beta \in \left] \alpha, \frac{1 + \alpha}{2} \right[$.
\end{proposition}

The function $f_0$ here satisfies the conditions of \cite[Theorem 1.5]{MR4687278}, yet its behavior near $z_c$ is very different from the example given above from \cite{ghoul2024nonlinear}. First of all, the scaling is different,
and the limit profile is still discontinuous here, while it was continuous in the example of \cite[Theorem 1.4]{ghoul2024nonlinear}. 

In  \cref{prop:ab}, the discontinuity is shown only for almost all $\beta \in
\left] \alpha, \frac{1 + \alpha}{2} \right[$, but it is likely that this holds for all $\beta$ in this range, although the proof seems technical.

If we add other terms of the form $\frac{\kappa}{| x |^{\delta}}$
in the initial data, we can nest these discontinuities. This is summarized in the following result.

\begin{theorem}
  \label{th:main} 
  There exist a strictly increasing sequence $(\alpha_n)_{n \in
  \mathbb{N}}$ with $0 < \alpha_n < 1$, and $\alpha_{\infty} = \lim_{n \to + \infty} \alpha_n <
  1$ such that, the solution of the viscous Burgers equation \cref{eq:burgers} with an initial data $f_0$ satisfying
   \begin{align}\label{eq:ic2}
    \begin{aligned}
 &  f_0 \in
C^1_{\tmop{loc}} (\mathbb{R}; \mathbb{R}), \\
& f_0 (x) = \sum_{n \in \mathbb{N}} \frac{2^{- n}}{| x |^{\alpha_n}} + o
     \left( \frac{1}{| x |} \right) \quad \text{as $x \to \pm \infty$}, \\
& f_0' (x) = \sum_{n \in \mathbb{N}}
     \frac{- \alpha_n 2^{- n}}{| x |^{1 + \alpha_n}} + o \left( \frac{1}{| x
     |^2} \right)  \quad \text{as $x \to \pm \infty$},
     \end{aligned}
    \end{align}
has the following asymptotic development: there exists two constants $y_+ \neq y_-$ such that, for any $N \in \mathbb N$ and $x \neq 0$, 
  \[ t^{\frac{\alpha_0}{1 + \alpha_0}} f \left( t^{\frac{\alpha_0}{1 +
     \alpha_0}} \left( z_c + t^{- \frac{\alpha_n - \alpha_0}{1 + \alpha_0}} x
     \right),t \right) - \sum_{n \leqslant N} t^{- \frac{\alpha_n - \alpha_0}{1
     + \alpha_0}} 2^{- n} Y_n (x) = o_{t \rightarrow + \infty} \left( t^{-
     \frac{\alpha_N - \alpha_0}{1 + \alpha_0}} \right), \]
where 
  \[ Y_n (x) \coloneqq \begin{cases}
       y_+^{\alpha_n} &\tmop{if } x > 0\\
       y_-^{\alpha_n} &\tmop{if } x < 0,
     \end{cases}
    \]
  for any $N \in \mathbb{N}$ and $x \in \mathbb{R}^{\ast}$,
\end{theorem}

This results in an example of an initial data for which the asymptotic profile
has a discontinuity, and although we can keep zooming on it, with new scales,
the discontinuity does not disappear after any arbitrary number of rescaling.

The rest of the paper is devoted to the proofs of these two results. We focus
first on the proof of  \cref{prop:ab}, and then we will show how to
adapt these arguments to show \cref{th:main}.

\section{Proof of \texorpdfstring{ \cref{prop:ab}}{Proposition 1.1}}

We define, as in {\cite[Eq. (2.2)]{MR4687278}}, the function
\[ H_t (y, z) \coloneqq \frac{- (z - y)^2}{4} - \frac{t^{- \frac{1 - \alpha}{1 +
   \alpha}}}{2} \int_0^{y t^{\frac{1}{1 + \alpha}}} f_0(\xi) \, \mathrm d \xi . \]
By the Hopf--Cole formula (see {\cite{MR42889,MR47234}} and also \cite[Section 4.4.1]{MR2597943}) and a change of variables, we have that the solution of \cref{eq:burgers} is given by 
\[ t^{\frac{\alpha}{1 + \alpha}} f \left( z t^{\frac{1}{1 + \alpha}}, t
   \right) = \frac{\displaystyle \int_{\mathbb{R}} t^{\frac{\alpha}{1 + \alpha}} f_0 \left(
   y t^{\frac{1}{1 + \alpha}} \right) e^{t^{\frac{1 - \alpha}{1 + \alpha}} H_t
   (y, z)} \, \d t}{\displaystyle \int_{\mathbb{R}} e^{t^{\frac{1 - \alpha}{1 + \alpha}} H_t (y,
   z)} \, \d t} . \]
By  \cite[Section 3.2.5]{MR4687278}, if $H_t (\cdot, z)$ has
a unique maximum, say $y_{\ast} (z, t)$, then
\[ t^{\frac{\alpha}{1 + \alpha}} f \left( z t^{\frac{1}{1 + \alpha}}, t
   \right) - t^{\frac{\alpha}{1 + \alpha}} f_0 \left( t^{\frac{1}{1 + \alpha}}
   y_{\ast} (z, t) \right) \longrightarrow 0 \quad \text{ as $t \nearrow +\infty$}.\]

In {\cite[Section 3.2.4]{MR4687278}}, it was shown that $H_t (\cdot, z)$ has
indeed one and only one maximum for all $z \neq z_c$ and $t$ large enough (depending on $| z - z_c |$). 

Here, we want to show that this is still the case if $z =
z_c + t^{- \frac{\beta - \alpha}{1 + \alpha}} x$, with $x \neq 0$ fixed and $t \nearrow + \infty$. To this end, we compute
\begin{equation}
  2 \partial_y H_t (y, z) = z - y - t^{\frac{\alpha}{1 + \alpha}} f_0 \left( y
  t^{\frac{1}{1 + \alpha}} \right) . \label{eq:ok}
\end{equation}
Owing to {\cite[Section 3.2.4]{MR4687278}}, there exists $a > 0$ such that, for
$z \in [z_c - a, z_c + a]$, among the solutions of $\partial_y H_t (y, z) = 0$
there are only two, denoted by $y_{\pm} (z, t)$, that are candidates to reach the
maximum of $H_t (\cdot, z)$ for $t$ large enough. Both are far away from $0$ and,
in the limit $t \rightarrow + \infty$, they converge to limits $y^{\ast}_{\pm}
(z)$ that are finite and solutions of the following implicit equation:\footnote{~Since $t^{\frac{\alpha}{1 + \alpha}} f_0 \left( y t^{\frac{1}{1 +
\alpha}} \right) \longrightarrow \frac{\kappa_1}{| y |^{\alpha}}$ when $t
\nearrow + \infty$ if $y \neq 0$, this equation comes from the fact that the right hand side of \cref{eq:ok} vanishes at a maximum.}
\[ z = y^{\ast}_{\pm} (z) + \frac{\kappa_1}{| y^{\ast}_{\pm} (z) |^{\alpha}}. 
\]

For $z = z_c$, we have 
\[\lim_{t \rightarrow + \infty} H_t (y_+ (z_c, t), z_c) -
H_t (y_- (z_c, t), z_c) = 0, 
\]
but, for $z < z_c$, this limit is strictly negative and, for $z > z_c$, it is strictly positive.

With the additional structure of $f_0$ prescribed in \cref{eq:ic1}, we can refine the estimate of $y_{\pm} (z, t)$ when $t$ is large.

\begin{lemma}
  \label{lm:21} 
  If $f_0$ satisfies \cref{eq:ic1}, we have
  that
  \begin{align}\label{eq:ypm} y_{\pm} (z, t) = y^{\ast}_{\pm} (z) + t^{- \frac{\beta - \alpha}{1 +
     \alpha}} \frac{\kappa_2}{| y^{\ast}_{\pm} (z) |^{\beta} \left(
     \frac{\kappa_1 \alpha}{| y^{\ast}_{\pm} (z) |^{1 + \alpha}} - 1 \right)}
     + o_{t \rightarrow + \infty} \left( t^{- \frac{\beta - \alpha}{1 +
     \alpha}} \right) . \end{align}
\end{lemma}

\begin{remark}
From {\cite[Section 3.1]{MR4687278}}, we have that ${y_{\pm}^{\ast}}' (z) \neq 0$, and differentiating $z = y^{\ast}_{\pm} (z) + \frac{\kappa_1}{| y^{\ast}_{\pm} (z)
|^{\alpha}}$ with respect to $z$ leads to $\frac{\kappa_1 \alpha}{|
y^{\ast}_{\pm} (z) |^{1 + \alpha}} - 1 = \frac{- 1}{{y_{\pm}^{\ast}}' (z)}
\neq 0$, hence the denominator in \cref{eq:ypm} never vanishes.
\end{remark}

\begin{proof}
  Since $y_{\pm} (z, t)$ is far away from $0$ for $t$ large
  enough (owing to {\cite[Section 3.2.3]{MR4687278}}), we have
  \[ t^{\frac{\alpha}{1 + \alpha}} f_0 \left( y_{\pm} (z, t) t^{\frac{1}{1 +
     \alpha}} \right) = \frac{\kappa_1}{| y_{\pm} (z, t) |^{\alpha}} +
     \frac{\kappa_2 t^{- \frac{\alpha - \beta}{1 + \alpha}}}{| y_{\pm} (z, t)
     |^{\beta}} + o_{t \rightarrow + \infty} \left( t^{- \frac{\alpha -
     \beta}{1 + \alpha}} \right) . \]
  Since $\partial_y H_t (y_{\pm} (z, t), z) = 0$ and $z = y^{\ast}_{\pm}
  (z) + \frac{\kappa_1}{| y^{\ast}_{\pm} (z) |^{\alpha}}$,  we compute that
  \[ y^{\ast}_{\pm} (z) - y_{\pm} (z, t) + \frac{\kappa_1}{| y^{\ast}_{\pm}
     (z) |^{\alpha}} - \frac{\kappa_1}{| y_{\pm} (z, t) |^{\alpha}} -
     \frac{\kappa_2 t^{- \frac{\alpha - \beta}{1 + \alpha}}}{| y_{\pm} (z, t)
     |^{\beta}} = o_{t \rightarrow + \infty} \left( t^{- \frac{\alpha -
     \beta}{1 + \alpha}} \right) . \]
  Observing that
  \[ \frac{\kappa_1}{| y^{\ast}_{\pm} (z) |^{\alpha}} - \frac{\kappa_1}{|
     y_{\pm} (z, t) |^{\alpha}} = \frac{\alpha \kappa_1}{| y^{\ast}_{\pm} (z)
     |^{1 + \alpha}} (y_{\pm} (z, t) - y^{\ast}_{\pm} (z)) + O ((y_{\pm} (z,
     t) - y^{\ast}_{\pm} (z))^2), \]
we conclude the proof.
\end{proof}

We can now compute the first correction to the values of $H_t (y_{\pm} (z, t),
z)$ for $z$ close to $z_c$.

\begin{lemma}
  \label{lm:22} For any $x \in \mathbb{R}$,
  \begin{align*}
      & H_t \left( y_{\pm} \left( z_c + t^{- \frac{\beta - \alpha}{1 +
    \alpha}} x, t \right), z_c + t^{- \frac{\beta - \alpha}{1 + \alpha}} x
    \right)\\
    & =  H_{\infty} (y^{\ast}_{\pm} (z_c)) + t^{- \frac{\beta - \alpha}{1 +
    \alpha}} x \left( \frac{z_c - y^{\ast}_{\pm} (z_c)}{2} \right) + o_{t
    \rightarrow + \infty} \left( t^{- \frac{\beta - \alpha}{1 + \alpha}}
    \right),
  \end{align*}
where $H_\infty \coloneqq  \lim_{t \rightarrow + \infty } H_t$.
\end{lemma}

\begin{proof}
  Since $\partial_y H_t (y_{\pm} (z, t), z) = 0$, we compute 
  \begin{align*}
    H_t (y_{\pm} (z, t), z)
    & =  H_{\infty} (y^{\ast}_{\pm} (z_c), z_c) + (z - z_c) \partial_z H_t
    (y_{\pm} (z, t), z)\\
    & \quad +  (y_{\pm} (z, t) - y^{\ast}_{\pm} (z_c))^2 \partial_y^2 H_t (y_{\pm}
    (z, t), z)\\
    & \quad +  O ((z - z_c)^2) + O ((y_{\pm} (z, t) - y^{\ast}_{\pm} (z_c))^2) .
  \end{align*}
  By \cref{lm:21}, we have $(y_{\pm} (z, t) - y^{\ast}_{\pm} (z_c))^2
  \partial_y^2 H_t (y_{\pm} (z, t), z) = o_{t \rightarrow + \infty} \left(
  t^{- \frac{\beta - \alpha}{1 + \alpha}} \right)$ and we conclude by taking
  $z = z_c + t^{- \frac{\beta - \alpha}{1 + \alpha}} x$ in $\partial_z H_t (y,
  z) = \frac{- (z - y)}{2}$.
\end{proof}

Since, owing to {\cite[Section 3.1]{MR4687278}}, $H_{\infty} (y^{\ast}_+ (z_c)) = H_{\infty} (y^{\ast}_- (z_c))$ and $z_c - y^{\ast}_+ (z_c) < 0$ while $z_c - y_-^{\ast} (z_c) > 0$, we infer the following corollary.

\begin{corollary}
  \label{cor:23}For any $\mu > 0$, there exists $K_{\nu} > 0$ such that, for $z
  = z_c + t^{- \frac{\beta - \alpha}{1 + \alpha}} x$ with $x > \mu$, the
  maximum of $H_t (\cdot, z)$ is reached, for $t$ large enough (depending on
  $\mu$) only at $y_+ (z, t)$, and outside of a neighborhood of size
  independent of time around $y_+ (z, t)$, the values of $H_t (\cdot, z)$ are away 
  from its maximum by a distance of at least $K_{\nu} t^{- \frac{\beta -
  \alpha}{1 + \alpha}}$.  The same holds for $x < - \mu$ by replacing $y_+$ by $y_-$.
\end{corollary}

Using these preliminary results, we can now show the following first order development of the solution to the viscous Burgers equation \cref{eq:burgers}.

\begin{lemma}
  \label{lm:24} For any $\mu > 0$, there exists $K_{\mu} > 0$ such that, for any
  $z \in \mathbb{R}$ with $| z - z_c | \geqslant \mu t^{- \frac{\beta -
  \alpha}{1 + \alpha}}$, if $z > z_c$ we have
  \[ \left| t^{\frac{\alpha}{1 + \alpha}} f \left( t^{\frac{\alpha}{1 +
     \alpha}} z \right) - t^{\frac{\alpha}{1 + \alpha}} f_0 \left( y_+ (z, t)
     t^{\frac{1}{1 + \alpha}} \right) \right| \leqslant K_{\mu} t^{- \frac{1 -
     \alpha}{2 (1 + \alpha)}}, \]
  and the same is true if $z < z_c$ up to replacing $y_+$ by $y_-$.
\end{lemma}

\begin{proof}
  By {\cite[Section 3.2.5]{MR4687278}}, 
  \[ \left| t^{\frac{\alpha}{1 + \alpha}} f \left( t^{\frac{\alpha}{1 +
     \alpha}} z \right) - t^{\frac{\alpha}{1 + \alpha}} f_0 \left( y (z, t)
     t^{\frac{1}{1 + \alpha}} \right) \right| \leqslant K \left( t^{-
     \frac{1}{2} + \frac{\alpha}{1 + \alpha}} + e^{- \frac{\nu}{2} t^{\frac{1
     - \alpha}{1 + \alpha}}} \right) \]
  provided that $y (z, t)$ is the only maximum of $H_t (\cdot, z)$, and where
  $\nu$ is the distance of $\max H_t \left(\cdot, z_c + t^{- \frac{\beta -
  \alpha}{1 + \alpha}} x \right)$ to any other value of $H_t (\cdot, z)$ except in
  a neighborhood of $y (z, t)$. In {\cite[Section 3.2.5]{MR4687278}}, this was done with $\nu$  uniform in time, but the proof still works with $\nu = \mu t^{- \frac{\beta
  - \alpha}{1 + \alpha}}$, as it is the case here due to \cref{cor:23}. Moreover,  since
  \[ \nu t^{\frac{1 - \alpha}{1 + \alpha}} = \mu t^{\frac{1 - \beta}{1 +
     \alpha}}, \quad \text{  with $1 - \beta > 0$,} \]
 we have
  \[ K \left( t^{- \frac{1}{2} + \frac{\alpha}{1 + \alpha}} + e^{-
     \frac{\nu}{2} t^{\frac{1 - \alpha}{1 + \alpha}}} \right) \leqslant
     K_{\mu} t^{\frac{\alpha - 1}{2 (1 + \alpha)}} \]
  where $K_{\mu} > 0$ depends on $\mu$.
\end{proof}

To establish \cref{prop:ab}, we need to compute the development of
$t^{\frac{\alpha}{1 + \alpha}} f \left( t^{\frac{\alpha}{1 + \alpha}} z
\right)$ up to a $o_{t \rightarrow + \infty} \left( t^{- \frac{\beta -
\alpha}{1 + \alpha}} \right)$. For the approximation of \cref{lm:24}
to be good enough to do so, we require that $\frac{\beta - \alpha}{1 + \alpha} < \frac{1 - \alpha}{2 (1 + \alpha)}$, that is $\beta < \frac{1 + \alpha}{2}$, which is an hypothesis of \cref{prop:ab}\footnote{~Below this threshold, we believe that there is another correction, independent of $\kappa_2, \beta$,
that needs to be considered. In a similar fashion, for the heat equation with integrable initial data, the first order in the long-time asymptotics is of size
$\frac{1}{\sqrt{t}}$ and the next term is of size $\frac{1}{t}$.}. 

\begin{proof}[Proof of \cref{prop:ab}]
  Recalling \cref{lm:24}, we have to develop
  \[ t^{\frac{\alpha}{1 + \alpha}} f_0 (y_{\pm} (z, t)) = \frac{\kappa_1}{|
     y_{\pm} (z, t) |^{\alpha}} + \frac{\kappa_2 t^{- \frac{\beta - \alpha}{1
     + \alpha}}}{| y_{\pm} (z, t) |^{\beta}} + o_{t \rightarrow + \infty}
     \left( t^{- \frac{\beta - \alpha}{1 + \alpha}} \right). \]

On the other hand, by \cref{lm:21}, 
  \[ y_{\pm} (z, t) = y^{\ast}_{\pm} (z) + t^{- \frac{\beta - \alpha}{1 +
     \alpha}} \frac{\kappa_2}{| y^{\ast}_{\pm} (z) |^{\beta} \left(
     \frac{\kappa_1 \alpha}{| y^{\ast}_{\pm} (z) |^{1 + \alpha}} - 1 \right)}
     + o_{t \rightarrow + \infty} \left( t^{- \frac{\beta - \alpha}{1 +
     \alpha}} \right) . \]
     
  By {\cite[Theorem 1.5]{MR4687278}}, $\mathfrak{p} (z) = \frac{\kappa_1}{|
  y^{\ast}_{\pm} (z) |^{\alpha}}$. Hence, we compute 
  \begin{align*}
    &   t^{\frac{\alpha}{1 + \alpha}} f_0 (y_{\pm} (z, t)) -\mathfrak{p}
    (z)\\
    & =  \frac{\kappa_1}{| y_{\pm} (z, t) |^{\alpha}} - \frac{\kappa_1}{|
    y^{\ast}_{\pm} (z) |^{\alpha}} + \frac{\kappa_2 t^{- \frac{\beta -
    \alpha}{1 + \alpha}}}{| y^{\ast}_{\pm} (z) |^{\beta}} + o_{t \rightarrow +
    \infty} \left( t^{- \frac{\beta - \alpha}{1 + \alpha}} \right)\\
    & =  t^{- \frac{\beta - \alpha}{1 + \alpha}} \frac{\kappa_2}{|
    y^{\ast}_{\pm} (z) |^{\beta}} \left( 1 + \frac{\alpha \kappa_1}{\alpha
    \kappa_1 - | y^{\ast}_{\pm} (z) |^{1 + \alpha}} \right) + o_{t \rightarrow
    + \infty} \left( t^{- \frac{\beta - \alpha}{1 + \alpha}} \right) .
  \end{align*}
  Since $y^{\ast}_{\pm} \left( z_c + t^{- \frac{\beta - \alpha}{1 + \alpha}} x
  \right) = y^{\ast}_{\pm} (z_c) + O_{t \rightarrow + \infty} \left( t^{-
  \frac{\beta - \alpha}{1 + \alpha}} \right)$, this proves that
  \begin{align*}
    &   t^{\frac{\beta - \alpha}{1 + \alpha}} \left( t^{\frac{\alpha}{1 +
    \alpha}} f_0 \left( y_{\pm} \left( z_c + t^{- \frac{\beta - \alpha}{1 +
    \alpha}} x, t \right) \right) -\mathfrak{p} \left( z_c + t^{- \frac{\beta
    - \alpha}{1 + \alpha}} x \right) \right)\\
    & \longrightarrow  \frac{\kappa_2}{| y^{\ast}_{\pm} (z_c) |^{\beta}} \left(
    1 + \frac{\alpha \kappa_1}{\alpha \kappa_1 - | y^{\ast}_{\pm} (z_c) |^{1 +
    \alpha}} \right) \eqqcolon \mathcal{P}_{\pm} (\beta) \quad \text{as $t \nearrow + \infty$},
  \end{align*}
  for $x \neq 0$, taking $y_+^{\ast}$ if $x > 0$
  and $y_-^{\ast}$ if $x < 0$.
  
  Since $z = y^{\ast}_{\pm} (z) + \frac{\kappa_1}{| y^{\ast}_{\pm} (z)
  |^{\alpha}}$, we have $\frac{{y_{\pm}^{\ast}}^{1 + \alpha}
  (z)}{{y_{\pm}^{\ast}}' (z)} {= y_{\pm}^{\ast}}^{1 + \alpha} (z) - \alpha
  \kappa_1$ and $\frac{\alpha \kappa_1 {y_{\pm}^{\ast}}'
  (z)}{{y_{\pm}^{\ast}}^{1 + \alpha} (z)} {= - 1 + y_{\pm}^{\ast}}' (z)$, which leads to 
  \[ \mathcal{P}_{\pm} (\beta) = \frac{\kappa_2}{| y^{\ast}_{\pm} (z_c)
     |^{\beta}} \left( {2 - y^{\ast}_{\pm}}' (z_c) \right) . \]

  By {\cite[Section 3.1]{MR4687278}}, we have ${y^{\ast}_-}' (z_c) < 0$ thus $\mathcal{P}_-
  (\beta) \neq 0$. We compute that
  \[ \mathcal{P}_{\pm}' (\beta) = - \ln (y^{\ast}_{\pm} (z_c))
     \mathcal{P}_{\pm} (\beta), \]
  hence, in particular if $\mathcal{P}_+ (\beta_0) =\mathcal{P}_- (\beta_0)
  \neq 0$, then since $y_+^{\ast} (z_c) \neq y_-^{\ast} (z_c)$ we have
  $\mathcal{P}_+' (\beta_0) \neq \mathcal{P}_-' (\beta_0)$ and thus there
  exists $\beta$ in a neighborhood of $\beta_0$ such that $\mathcal{P}_+
  (\beta_0) \neq \mathcal{P}_- (\beta_0) .$ This implies that the $\beta \in
  \left] \alpha, \frac{1 + \alpha}{2} \right[$ such that $\mathcal{P}_+
  (\beta) \neq \mathcal{P}_- (\beta)$ are dense. For these values, the profile
  is discontinuous at $x = 0$.
\end{proof}

\section{Proof of \texorpdfstring{\cref{th:main}}{Theorem 1.2}}

Following the same steps as in the proof of  \cref{prop:ab}, we can prove \cref{th:main} as well.

\begin{proof}[Proof of \cref{th:main}]
  First, we note that there exists an increasing sequence $(\alpha_n)_{n \in \mathbb{N}}$ such that $0 <
  \alpha_0 < 1$, $\alpha_{\infty} = \lim_{n \rightarrow \infty} \alpha_n < 1$,   and
  \begin{align}\label{ass:alpha}\alpha_1 > \frac{\alpha_{\infty} +
  \alpha_0}{2}, \qquad \alpha_n < \frac{1 + \alpha_0}{2}, \quad \text{for all $n \in \mathbb{N}$.} 
  \end{align}

Then we compute 
    \begin{equation}
    t^{\frac{\alpha_0}{1 + \alpha_0}} f_0 \left( y_{\pm} (z, t) t^{\frac{1}{1
    + \alpha_0}} \right) = \sum_{n \in \mathbb{N}^{\ast}} \frac{2^{- n} t^{-
    \frac{\alpha_n - \alpha_0}{1 + \alpha_0}}}{| y_{\pm} (z, t) |^{\alpha_n}}
    + o_{t \rightarrow + \infty} \left( t^{- \frac{1 - \alpha_0}{1 +
    \alpha_0}} \right). \label{eq:21}
  \end{equation}
  Since $0 = 2 \partial_y H_t (y_{\pm} (z, t), z) = z - y_{\pm} (z, t) -
  t^{\frac{\alpha}{1 + \alpha}} f_0 \left( y_{\pm} (z, t) t^{\frac{1}{1 +
  \alpha}} \right)$, we deduce, similarly as in \cref{lm:21}, that 
  \begin{equation}
    y_{\pm} (z, t) = y^{\ast}_{\pm} (z) + \sum_{n \in \mathbb{N}^{\ast}} t^{-
    \frac{\alpha_n - \alpha_0}{1 + \alpha_0}} Y_{\pm, n} (z) + o_{t
    \rightarrow + \infty} \left( t^{- \frac{\alpha_{\infty} - \alpha_0}{1 +
    \alpha_0}} \right), \label{eq:22}
  \end{equation}
  where 
    \[ Y_{\pm, n} (z) \coloneqq  \frac{2^{- n}}{| y^{\ast}_{\pm} (z) |^{\beta} \left(
     \frac{\alpha_0}{| y^{\ast}_{\pm} (z) |^{1 + \alpha_0}} - 1 \right)}. \] 
     
  In obtaining \cref{eq:22}, we used that   
  \[\frac{\alpha_n - \alpha_0}{1 + \alpha_0} + \frac{\alpha_m - \alpha_0}{1 +
  \alpha_0} > \frac{\alpha_{\infty} - \alpha_0}{1 + \alpha_0} \quad \text{ for all $(n, m)
  {\in \mathbb{N}^{\ast}}^2$, \ $n \neq m$.}
  \]
  Since the $(\alpha_n)_{n \in
  \mathbb{N}}$ are increasing, this holds, in particular, if $\alpha_1 >  \frac{\alpha_{\infty} + \alpha_0}{2}$, which we imposed.
  
  As in the proof of \cref{lm:22}, we check that
  \begin{align*}
    &  H_t \left( y_{\pm} \left( z_c + t^{- \frac{\alpha_n - \alpha_0}{1 +
    \alpha_0}} x, t \right), z_c + t^{- \frac{\alpha_n - \alpha_0}{1 +
    \alpha_0}} x \right)\\
    & =  H_{\infty} (y^{\ast}_{\pm} (z_c)) + t^{- \frac{\alpha_n -
    \alpha_0}{1 + \alpha_0}} x \left( \frac{z_c - y^{\ast}_{\pm} (z_c)}{2}
    \right) + o_{t \rightarrow + \infty} \left( t^{- \frac{\alpha_n -
    \alpha_0}{1 + \alpha_0}} \right),
  \end{align*}
  provided that $\alpha_1 > \frac{\alpha_{\infty} + \alpha_0}{2}$.
  
  We deduce  that, for $\mu > 0$, the maximum is reached at 
  \begin{align*} 
  \displaystyle y_+ \left( z_c + t^{-
  \frac{\alpha_n - \alpha_0}{1 + \alpha_0}} x, t \right)  &\qquad\text{ if $x > \mu$}, \\ 
  \displaystyle y_- \left( z_c + t^{- \frac{\alpha_n - \alpha_0}{1 + \alpha_0}} x, t
  \right) &\qquad\text{ if $x < - \mu$}
  \end{align*} 
  if $t$ is large enough (depending on $\mu$ and $n$),
  and the difference between the two maxima is of size $ K_{\mu} t^{-
  \frac{\alpha_n - \alpha_0}{1 + \alpha_0}},$  for some $K_{\mu} > 0$ (depending 
  on $\mu$ and $n$).
  
  Therefore, the proof of \cref{lm:24} still holds and we have that, if $z > z_c$, 
  \[ \left| t^{\frac{\alpha_0}{1 + \alpha_0}} f \left( t^{\frac{\alpha_0}{1 +
     \alpha_0}} z \right) - t^{\frac{\alpha_0}{1 + \alpha_0}} f_0 \left( y_+
     (z, t) t^{\frac{1}{1 + \alpha_0}} \right) \right| \leqslant K_{\mu} t^{-
     \frac{1 - \alpha_0}{2 (1 + \alpha_0)}} \leqslant K_{\mu} t^{-
     \frac{\alpha_{\infty} - \alpha_0}{1 + \alpha_0}} \quad \text{for $t \ge 1$,} \]
   (and the same is true if $z < z_c$ up to replacing $y_+$ by $y_-$),
  since $\alpha_{\infty} \leqslant \frac{1 + \alpha_0}{2}$.
  
  Combining this with \cref{eq:21} and \cref{eq:22}, since   $\alpha_1 > \frac{\alpha_{\infty} + \alpha_0}{2}$, we conclude that 
  \[ t^{\frac{\alpha_0}{1 + \alpha_0}} f \left( t^{\frac{\alpha_0}{1 +
     \alpha_0}} \left( z_c + t^{- \frac{\alpha_n - \alpha_0}{1 + \alpha_0}} x
     \right) \right) - \sum_{n \leqslant N} \frac{2^{- n} t^{- \frac{\alpha_n
     - \alpha_0}{1 + \alpha_0}}}{| y^{\ast}_{\pm} (z_c) |^{\alpha_n}} = o_{t
     \rightarrow + \infty} \left( t^{- \frac{\alpha_N - \alpha_0}{1 +
     \alpha_0}} \right) \]
  for any $N \in \mathbb{N}$ and $x \in \mathbb{R}^{\ast}$.
\end{proof}

\section*{Acknowledgments}

N.~De Nitti is a member of the Gruppo Nazionale per l'Analisi Matematica, la Probabilità e le loro Applicazioni (GNAMPA) of the Istituto Nazionale di Alta Matematica (INdAM). He has been funded by the Swiss State Secretariat for Education, Research and Innovation (SERI) under contract number MB22.00034 through the project TENSE. He also acknowledges the kind hospitality of CY Cergy Paris University, where part of this work was carried out. 

\vspace{0.5cm}

\printbibliography
\vfill 
\end{document}